\documentclass[11pt,twoside]{amsart}
\usepackage{graphicx}
\usepackage{amssymb}
\usepackage{amsmath}
\usepackage{amsthm}
\usepackage{amsfonts}
\usepackage{amsmath, amsfonts, amssymb}
\newtheorem{theorem}{Theorem}[section]

\newtheorem{remark}[theorem]{Remark}

\newtheorem{definition}[theorem]{Definition}
\numberwithin{equation}{section}

\begin{document}
	\begin{center}\small{In the name of Allah, the Beneficent, the Merciful.}\end{center}
	\vspace{1cm}
\title{The automorphism groups and derivation algebras of two-dimensional algebras}

\author{H.Ahmed$^1$, U.Bekbaev$^2$, I.Rakhimov$^3$}

\thanks{{\scriptsize
emails: $^1$houida\_m7@yahoo.com; $^2$bekbaev@iium.edu.my; $^3$rakhimov@upm.edu.my.}}
\maketitle
\begin{center}
\address{$^{1}$Department of Math., Faculty of Science, UPM, Selangor, Malaysia $\&$ Depart. of Math., Faculty of Science, Taiz University, Taiz, Yemen}\\
  \address{$^2$Department of Science in Engineering, Faculty of Engineering, IIUM, Kuala Lumpur, Malaysia}\\
\address{$^3$Department of Math., Faculty of Science $\&$ Institute for Mathematical Research (INSPEM), UPM, Serdang, Selangor, Malaysia}
\end{center}
\begin{abstract} The automorphisms groups and derivation algebras of all two-dimensional algebras over algebraically closed fields are described.
\end{abstract}
\section{Introduction}

In \cite{A} the authors have presented a complete list of isomorphism classes of two-dimensional algebras over algebraically closed fields, providing a list of canonical representatives of their structure constant's matrices. In the present paper we describe the groups of automorphisms and derivation algebras of all those listed algebras.

In fact, the automorphism groups of all $2$-dimensional algebras have been given earlier in \cite{I} up to isomorphism. In contrast to that we provide an explicit realization of the automorphism groups for the all listed canonical algebras (see \cite{A}).

The first part of this paper (Sections 2 and 3) is devoted to the description of the groups of automorphisms and the second part (Section 4) deals with the derivation algebras. In each case we consider problems over algebraically closed fields of characteristic not $2,3$, characteristic $2$ and characteristic $3$ separately according to classification results of \cite{A}.

\section{ Preliminaries}

Let $\mathbb{F}$ be any field, for matrices $A=(a_{ij})$, $B$ over $\mathbb{F}$, as usual $A\otimes B$ stand for the block-matrix with blocks $(a_{ij}B).$

\begin{definition} A vector space  $\mathbb{A}$ over $\mathbb{F}$ with multiplication $\cdot :\mathbb{A}\otimes \mathbb{A}\rightarrow \mathbb{A}$ given by $(\mathbf{u},\mathbf{v})\mapsto \mathbf{u}\cdot \mathbf{v}$ such that \[(\alpha\mathbf{u}+\beta\mathbf{v})\cdot \mathbf{w}=\alpha(\mathbf{u}\cdot \mathbf{w})+\beta(\mathbf{v}\cdot \mathbf{w}),\ \ \mathbf{w}\cdot (\alpha\mathbf{u}+\beta\mathbf{v})=\alpha(\mathbf{w}\cdot \mathbf{u})+\beta(\mathbf{w}\cdot \mathbf{v})\] whenever $\mathbf{u}, \mathbf{v}, \mathbf{w}\in \mathbb{A}$ and $\alpha, \beta\in \mathbb{F}$, is said to be an algebra.\end{definition}

\begin{definition} Two algebras $\mathbb{A}$ and $\mathbb{B}$ are called isomorphic if there is an invertible linear map  $\mathbf{f}:\mathbb{A}\rightarrow \mathbb{B} $ such that \[\mathbf{f}(\mathbf{u}\cdot_{\mathbb{A}} \mathbf{v})=\mathbf{f}(\mathbf{u})\cdot_{\mathbb{B}} \mathbf{f}(\mathbf{v})\] whenever $\mathbf{u}, \mathbf{v}\in \mathbb{A}$ and $\alpha, \beta\in \mathbb{F}$.\end{definition}

\begin{definition} An invertible linear map  $\mathbf{f}:\mathbb{A}\rightarrow \mathbb{A} $ is said to be an automorphism if \[\mathbf{f}(\mathbf{u}\cdot \mathbf{v})=\mathbf{f}(\mathbf{u})\cdot \mathbf{f}(\mathbf{v})\] whenever $\mathbf{u}, \mathbf{v}\in \mathbb{A}$ and $\alpha, \beta\in \mathbb{F}$.\end{definition}

\begin{definition} A linear map  $\mathbf{d}:\mathbb{A}\rightarrow \mathbb{A} $ is said to be a derivation if \[\mathbf{d}(\mathbf{u}\cdot \mathbf{v})=\mathbf{d}(\mathbf{u})\cdot \mathbf{v}+\mathbf{u}\cdot \mathbf{d}(\mathbf{v})\] whenever $\mathbf{u}, \mathbf{v}\in \mathbb{A}$ and $\alpha, \beta\in \mathbb{F}$.\end{definition}

Let $\mathbb{A}$ be $m$-dimensional algebra over $\mathbb{F}$ and $e=(e^1,e^2,...,e^m)$ its basis. Then the bilinear map $\cdot$ is represented by a matrix $A=(A^k_{ij})\in M(m\times m^2;\mathbb{F})$ as follows \[\mathbf{u}\cdot \mathbf{v}=eA(u\otimes v),\] for $\mathbf{u}=eu,\mathbf{v}=ev,$
where $u = (u_1, u_2, ..., u_m)^T,$ and $v = (v_1, v_2, ..., v_m)^T$ are column coordinate vectors of $\mathbf{u}$ and $\mathbf{v},$ respectively.
The matrix $A\in M(m\times m^2;\mathbb{F})$ defined above is called the matrix of structural constants (MSC) of $\mathbb{A}$ with respect to the basis $e$. Further we assume that a basis $e$ is fixed and we do not make a difference between the algebra
$\mathbb{A}$ and its MSC $A$.

An automorphism $\mathbf{g}:\mathbb{A}\rightarrow \mathbb{A}$ as an invertible linear map is represented by an invertible matrix $g\in GL(m;\mathbb{F})$: $\mathbf{g}(\mathbf{u})=\mathbf{g}(eu)=egu$. Due to  \[\mathbf{g}(\mathbf{u}\cdot \mathbf{v})=\mathbf{g}(eA(u\otimes v)=eg(A(u\otimes v))=e(gA)(u\otimes v),\] and \[\mathbf{g}(\mathbf{u})\cdot \mathbf{g}(\mathbf{v})=(egu)\cdot (egv)=eA(gu\otimes gv)=eAg^{\otimes 2}(u\otimes v)\] the property $\mathbf{g}(\mathbf{u}\cdot \mathbf{v})=\mathbf{g}(\mathbf{u})\cdot \mathbf{g}(\mathbf{v})$ is equivalent to \begin{equation}\label{1}gA=Ag^{\otimes 2}.\end{equation}

An derivation $\mathbf{d}:\mathbb{A}\rightarrow \mathbb{A}$ as a linear map is represented by a matrix $d\in M(m;\mathbb{F})$ as follows $\mathbf{d}(\mathbf{u})=\mathbf{d}(eu)=edu$. Due to  \[\mathbf{d}(\mathbf{u}\cdot \mathbf{v})=\mathbf{d}(eA(u\otimes v))=ed(A(u\otimes v))=e(dA)(u\otimes v),\] and \[\mathbf{d}(\mathbf{u})\cdot \mathbf{v}+\mathbf{u}\cdot \mathbf{d}(\mathbf{v})=(edu)\cdot (ev)+(eu)\cdot (edv)=eA(du\otimes v)+eA(u\otimes dv)=\]
\[e(A(d\otimes I)(u\otimes v)+A(I\otimes d)(u\otimes v))= eA((d\otimes I)+(I\otimes d))(u\otimes v)\]  the property $\mathbf{d}(\mathbf{u}\cdot \mathbf{v})=\mathbf{d}(\mathbf{u})\cdot \mathbf{v}+\mathbf{u}\cdot \mathbf{d}(\mathbf{v})$ is equivalent to \begin{equation}\label{2}dA=A(d\otimes I+I\otimes d),\end{equation}
where $I$ stands for the identity matrix.

If $e'=(e'^1,e'^2,...,e'^m)$ is another basis of $\mathbb{A}$, $e'g=e$ with $g\in G=GL(m;\mathbb{F})$, and  $A'$ is MSC of $\mathbb{A}$ with respect to $e'$ then it is known that
\begin{equation}\label{3}A'=gA(g^{-1})^{\otimes 2}\end{equation} is valid. Thus, the isomorphism of algebras $\mathbb{A}$ and $\mathbb{B}$ over $\mathbb{F}$ given above  now can be rewritten as follows.

\begin{definition} Two $m$-dimensional algebras $\mathbb{A}$, $\mathbb{B}$ over $\mathbb{F}$, given by
	their matrices of structure constants $A$, $B$, are said to be isomorphic if $B=gA(g^{-1})^{\otimes 2}$ holds true for some $g\in GL(m;\mathbb{F})$.\end{definition}

Further we consider only the case $m=2$ and for the simplicity we use \[A=\left(\begin{array}{cccc} \alpha_1 & \alpha_2 & \alpha_3 &\alpha_4\\ \beta_1 & \beta_2 & \beta_3 &\beta_4\end{array}\right)\] for MSC, where
$\alpha_1, \alpha_2, \alpha_3, \alpha_4, \beta_1, \beta_2, \beta_3, \beta_4$ stand for any elements of $\mathbb{F}$.

Due to \cite{A} we have the following classification theorems according to $Char(\mathbb{F})\neq 2,3,$ $Char(\mathbb{F})=2$ and $Char(\mathbb{F})=3$ cases, respectively.
\begin{theorem}\label{thm1} Over an algebraically closed field $\mathbb{F}$ $(Char(\mathbb{F})\neq 2$ and $3)$, any non-trivial $2$-dimensional algebra is isomorphic to only one of the following algebras listed by their matrices of structure constants:
	\begin{itemize}
	\item $A_{1}(\mathbf{c})=\left(
	\begin{array}{cccc}
	\alpha_1 & \alpha_2 &\alpha_2+1 & \alpha_4 \\
	\beta_1 & -\alpha_1 & -\alpha_1+1 & -\alpha_2
	\end{array}\right),\ \mbox{where}\ \mathbf{c}=(\alpha_1, \alpha_2, \alpha_4, \beta_1)\in \mathbb{F}^4,$
	\item $A_{2}(\mathbf{c})=\left(
	\begin{array}{cccc}
	\alpha_1 & 0 & 0 & 1 \\
	\beta _1& \beta _2& 1-\alpha_1&0
	\end{array}\right)\simeq \left(
	\begin{array}{cccc}
	\alpha_1 & 0 & 0 & 1 \\
	-\beta _1& \beta _2& 1-\alpha_1&0
	\end{array}\right),\ \mbox{where}\ \mathbf{c}=(\alpha_1, \beta_1, \beta_2)\in \mathbb{F}^3,$
	\item $A_{3}(\mathbf{c})=\left(
	\begin{array}{cccc}
	0 & 1 & 1 & 0 \\
	\beta _1& \beta _2 & 1&-1
	\end{array}\right),\ \mbox{where}\ \mathbf{c}=(\beta_1, \beta_2)\in \mathbb{F}^2,$
	\item $A_{4}(\mathbf{c})=\left(
	\begin{array}{cccc}
	\alpha _1 & 0 & 0 & 0 \\
	0 & \beta _2& 1-\alpha _1&0
	\end{array}\right),\ \mbox{where}\ \mathbf{c}=(\alpha_1, \beta_2)\in \mathbb{F}^2,$
	\item $A_{5}(\mathbf{c})=\left(
	\begin{array}{cccc}
	\alpha_1& 0 & 0 & 0 \\
	1 & 2\alpha_1-1 & 1-\alpha_1&0
	\end{array}\right),\ \mbox{where}\ \mathbf{c}=\alpha_1\in \mathbb{F},$
	\item $A_{6}(\mathbf{c})=\left(
	\begin{array}{cccc}
	\alpha_1 & 0 & 0 & 1 \\
	\beta _1& 1-\alpha_1 & -\alpha_1&0
	\end{array}\right)\simeq \left(
	\begin{array}{cccc}
	\alpha_1 & 0 & 0 & 1 \\
	-\beta _1& 1-\alpha_1 & -\alpha_1&0
	\end{array}\right),\ \mbox{where}\ \mathbf{c}=(\alpha_1, \beta_1)\in \mathbb{F}^2,$
	\item $A_{7}(\mathbf{c})=\left(
	\begin{array}{cccc}
	0 & 1 & 1 & 0 \\
	\beta_1& 1& 0&-1
	\end{array}\right),\ \mbox{where}\ \mathbf{c}=\beta_1\in \mathbb{F},$
	\item $A_{8}(\mathbf{c})=\left(
	\begin{array}{cccc}
	\alpha_1 & 0 & 0 & 0 \\
	0 & 1-\alpha_1 & -\alpha_1&0
	\end{array}\right),\ \mbox{where}\ \mathbf{c}=\alpha_1\in \mathbb{F},$
	\item $A_{9}=\left(
	\begin{array}{cccc}
	\frac{1}{3}& 0 & 0 & 0 \\
	1 & \frac{2}{3} & -\frac{1}{3}&0
	\end{array}\right),$
	\item $A_{10}=\left(
	\begin{array}{cccc}
	0 & 1 & 1 & 0 \\
	0 &0& 0 &-1
	\end{array}
	\right),$
	\item $A_{11}=\left(
	\begin{array}{cccc}
	0 & 1 & 1 & 0 \\
	1 &0& 0 &-1
	\end{array}
	\right),$
	\item $A_{12}=\left(
	\begin{array}{cccc}
	0 & 0 & 0 & 0 \\
	1 &0&0 &0\end{array}
	\right).$
\end{itemize}\end{theorem}

\begin{theorem}\label{thm2} Over an algebraically closed field $\mathbb{F}$ $(Char(\mathbb{F})=2)$, any non-trivial $2$-dimensional algebra is isomorphic to only one of the following algebras listed by their matrices of structure constants:
\begin{itemize}
\item $A_{1,2}(\mathbf{c})=\left(
\begin{array}{cccc}
\alpha_1 & \alpha_2 &\alpha_2+1 & \alpha_4 \\
\beta_1 & -\alpha_1 & -\alpha_1+1 & -\alpha_2
\end{array}\right),\ \mbox{where}\ \mathbf{c}=(\alpha_1, \alpha_2, \alpha_4, \beta_1)\in \mathbb{F}^4,$
\item $A_{2,2}(\mathbf{c})=\left(
\begin{array}{cccc}
\alpha_1 & 0 & 0 & 1 \\
\beta _1& \beta_2 & 1-\alpha_1&0
\end{array}\right),\ \mbox{where}\ \mathbf{c}=(\alpha_1, \beta_1, \beta_2)\in \mathbb{F}^3,$
\item $A_{3,2}(\mathbf{c})=\left(
\begin{array}{cccc}
\alpha_1 & 1 & 1 & 0 \\
0& \beta_2 & 1-\alpha_1&1
\end{array}\right),\ \mbox{where}\ \mathbf{c}=(\alpha_1, \beta_2)\in \mathbb{F}^2,$
\item $A_{4,2}(\mathbf{c})=\left(
\begin{array}{cccc}
\alpha _1 & 0 & 0 & 0 \\
0 & \beta_2 & 1-\alpha _1&0
\end{array}\right),\ \mbox{where}\ \mathbf{c}=(\alpha_1,\beta_2)\in \mathbb{F}^2,$
\item $A_{5,2}(\mathbf{c})=\left(
\begin{array}{cccc}
\alpha_1 & 0 & 0 & 0 \\
1 & 1 & 1-\alpha_1&0
\end{array}\right),\ \mbox{where}\ \mathbf{c}=\alpha_1\in \mathbb{F},$
\item $A_{6,2}(\mathbf{c})=\left(
\begin{array}{cccc}
\alpha_1 & 0 & 0 & 1 \\
\beta _1& 1-\alpha_1 & -\alpha_1&0
\end{array}\right),\ \mbox{where}\ \mathbf{c}=(\alpha_1, \beta_1)\in \mathbb{F}^2,$
\item $A_{7,2}(\mathbf{c})=\left(
\begin{array}{cccc}
\alpha_1 & 1 & 1 & 0 \\
0& 1-\alpha_1& -\alpha_1&-1
\end{array}\right),\ \mbox{where}\ \mathbf{c}=\alpha_1\in \mathbb{F},$
\item $A_{8,2}(\mathbf{c})=\left(
\begin{array}{cccc}
\alpha_1 & 0 & 0 & 0 \\
0 & 1-\alpha_1 & -\alpha_1&0
\end{array}\right),\ \mbox{where}\ \mathbf{c}=\alpha_1\in \mathbb{F},$
\item $A_{9,2}=\left(
\begin{array}{cccc}
1 & 0 & 0 & 0 \\
1 & 0 & 1&0
\end{array}\right),$
\item $A_{10,2}=\left(
\begin{array}{cccc}
0 & 1 & 1 & 0 \\
0 &0& 0 &-1
\end{array}
\right),$
\item $A_{11,2}=\left(
\begin{array}{cccc}
1 & 1 & 1 & 0 \\
0 &-1& -1 &-1
\end{array}
\right),$
\item $A_{12,2}=\left(
\begin{array}{cccc}
0 & 0 & 0 & 0 \\
1 &0&0 &0\end{array}
\right).$
\end{itemize}
\end{theorem}

\begin{theorem}\label{thm3} Over an algebraically closed field $\mathbb{F}$ $(Char(\mathbb{F})=3)$, any non-trivial $2$-dimensional algebra is isomorphic to only one of the following algebras listed by their matrices of structure constant matrices:
\begin{itemize}
\item $A_{1,3}(\mathbf{c})=\left(
\begin{array}{cccc}
\alpha_1 & \alpha_2 &\alpha_2+1 & \alpha_4 \\
\beta_1 & -\alpha_1 & -\alpha_1+1 & -\alpha_2
\end{array}\right),\ \mbox{where}\ \mathbf{c}=(\alpha_1, \alpha_2, \alpha_4, \beta_1)\in \mathbb{F}^4,$
\item $A_{2,3}(\mathbf{c})=\left(
\begin{array}{cccc}
\alpha_1 & 0 & 0 & 1 \\
\beta _1& \beta _2& 1-\alpha_1&0
\end{array}\right)\simeq\left(
\begin{array}{cccc}
\alpha_1 & 0 & 0 & 1 \\
-\beta _1& \beta _2& 1-\alpha_1&0
\end{array}\right),\ \mbox{where}\ \mathbf{c}=(\alpha_1, \beta_1, \beta_2)\in \mathbb{F}^3,$
\item $A_{3,3}(\mathbf{c})=\left(
\begin{array}{cccc}
0 & 1 & 1 & 0 \\
\beta _1& \beta _2 & 1&-1
\end{array}\right),\ \mbox{where}\ \mathbf{c}=(\beta_1, \beta_2)\in \mathbb{F}^2,$
\item $A_{4,3}(\mathbf{c})=\left(
\begin{array}{cccc}
\alpha _1 & 0 & 0 & 0 \\
0 & \beta _2& 1-\alpha _1&0
\end{array}\right),\ \mbox{where}\ \mathbf{c}=(\alpha_1, \beta_2)\in \mathbb{F}^2,$
\item $A_{5,3}(\mathbf{c})=\left(
\begin{array}{cccc}
\alpha_1& 0 & 0 & 0 \\
1 & -1-\alpha_1 & 1-\alpha_1&0
\end{array}\right),\ \mbox{where}\ \mathbf{c}=\alpha_1\in \mathbb{F},$
\item $A_{6,3}(\mathbf{c})=\left(
\begin{array}{cccc}
\alpha_1 & 0 & 0 & 1 \\
\beta _1& 1-\alpha_1 & -\alpha_1&0
\end{array}\right)\simeq \left(
\begin{array}{cccc}
\alpha_1 & 0 & 0 & 1 \\
-\beta _1& 1-\alpha_1 & -\alpha_1&0
\end{array}\right),\ \mbox{where}\ \mathbf{c}=(\alpha_1, \beta_1)\in \mathbb{F}^2,$
\item $A_{7,3}(\mathbf{c})=\left(
\begin{array}{cccc}
0 & 1 & 1 & 0 \\
\beta_1& 1& 0&-1
\end{array}\right),\ \mbox{where}\ \mathbf{c}=\beta_1\in \mathbb{F},$
\item $A_{8,3}(\mathbf{c})=\left(
\begin{array}{cccc}
\alpha_1 & 0 & 0 & 0 \\
0 & 1-\alpha_1 & -\alpha_1&0
\end{array}\right),\ \mbox{where}\ \mathbf{c}=\alpha_1\in \mathbb{F},$
\item $A_{9,3}=\left(
\begin{array}{cccc}
0 & 1& 1& 0 \\
1 &0&0 &-1\end{array}
\right),$
\item $A_{10,3}=\left(
\begin{array}{cccc}
0 & 1 & 1 & 0 \\
0 &0&0 &-1\end{array}
\right),$
\item $A_{11,3}=\left(
\begin{array}{cccc}
1 & 0 & 0 & 0 \\
1 &-1&-1 &0\end{array}
\right),$
\item $A_{12,3}=\left(
\begin{array}{cccc}
0 &0 &0 & 0 \\
1 &0& 0 &0
\end{array}
\right).$
\end{itemize}\end{theorem}
\begin{remark} In \cite{A} the class $A_{3,2}(\mathbf{c})$ should be understood as it is in this paper, as far as there is a type-mistake in this case in \cite{A}.\end{remark}

\section{ The groups of automorphisms of $2$-dimensional algebras}
 Due to (\ref{1}) for the group of automorphisms of an algebra $\mathbb{A}$ given by MSC $A\in M(2\times 4;\mathbb{F})$ one has
\begin{equation}\label{aut}
  Aut(A)=\{g\in GL(2;\mathbb{F}): \  gA-A(g\otimes g)=0\}.
\end{equation}

Therefore further we look only for nonsingular solutions $g=\left(
\begin{array}{cc}
	a & b \\
	c & d \\
\end{array}
\right)$ of the equation $gA-A(g\otimes g)=0.$ We consider this equation for each MSC $A$ from the theorems stated above separately and to avoid reparations we discuss only a few cases in details, the others can be understood easily from the discussions made. We use $\Delta$ for $\det(g)=ad-bc$.

\begin{theorem}\label{thm4} The automorphism groups of the algebras listed in Theorem $\ref{thm1}$ are given as follows:
	\begin{itemize}
\item $Aut(A_1(\alpha_1,\alpha_2,\alpha_4,\beta_1))=\{I\},$

\item $Aut(A_2(\alpha_1,\beta_1,\beta_2))=\{I\}$ if $\ \beta_1\neq 0 ,$
 $Aut(A_2(\alpha_1,0,\beta_2))=\left\{I, \left(
\begin{array}{cc}
1 & 0 \\
0 & -1 \\
\end{array}
\right)\right\},$
\item $Aut(A_3(\beta_1,\beta_2))=\{I\},$
\item $Aut(A_4(\alpha_1,\beta_2))=\left\{ \left(
\begin{array}{cc}
1 & 0 \\
0 & d \\
\end{array}
\right):\ d\neq 0 \right\}$ if $\beta_2\neq 2\alpha_1-1 ,$
\item $Aut(A_4(\alpha_1,2\alpha_1-1 ))=\left\{ \left(
\begin{array}{cc}
1 & 0 \\
c & d \\
\end{array}
\right):c \in\mathbb{F}, d\neq 0 \right\},$\\
\item $Aut(A_5(\alpha_1))=\left\{ \left(
\begin{array}{cc}
1 & 0 \\
c & 1 \\
\end{array}
\right): c \in\mathbb{F} \right\},$\\
\item $Aut(A_6(\alpha_1,\beta_1))=\{I\}$ if $\beta_1\neq 0 ,$
 $Aut(A_6(\alpha_1,0))=\left\{I, \left(
\begin{array}{cc}
1 & 0 \\
0 & -1 \\
\end{array}
\right),\right\},$
\item $Aut(A_7(\beta_1))=\{I\}.$
\item $Aut(A_8(\alpha_1))=\left\{\left(
	\begin{array}{cc}
	1 & 0 \\
	0 & d \\
	\end{array}
	\right):\ d\neq 0 \right\}$ if $\alpha_1\neq \frac{1}{3} ,$
 $Aut(A_8(\frac{1}{3}))=\left\{\left(
	\begin{array}{cc}
	1 & 0 \\
	c & d \\
	\end{array}
	\right):c \in\mathbb{F}, d\neq 0 \right\},$
\item $Aut(A_9)=\left\{\left(
	\begin{array}{cc}
	1 & 0 \\
	c & 1 \\
	\end{array}
	\right): c\in \mathbb{F} \right\},$
\item $Aut(A_{10})=\left\{\left(
	\begin{array}{cc}
	a & 0 \\
	0 & 1 \\
	\end{array}
	\right):\ a\neq 0 \right\},$
\item $Aut(A_{11})=\left\{I,\left(
\begin{array}{cc}
- 1 & 0 \\
0 & 1 \\
\end{array}
\right),\left(
\begin{array}{cc}
\frac{1}{2} &  \frac{\sqrt{3}}{2} \\
 \frac{\sqrt{3}}{2} & -\frac{1}{2}  \\
\end{array}
\right),\left(
\begin{array}{cc}
\frac{1}{2} & - \frac{\sqrt{3}}{2} \\
- \frac{\sqrt{3}}{2} & -\frac{1}{2}  \\
\end{array}
\right),\left(
\begin{array}{cc}
-\frac{1}{2} &  \frac{\sqrt{3}}{2} \\
- \frac{\sqrt{3}}{2} & -\frac{1}{2}  \\
\end{array}
\right),\left(
\begin{array}{cc}
-\frac{1}{2} & - \frac{\sqrt{3}}{2} \\
 \frac{\sqrt{3}}{2} & -\frac{1}{2}  \\
\end{array}
\right)\right\},$
\item $Aut(A_{12})=\left\{\left(
	\begin{array}{cc}
	a & 0 \\
	c & a^2 \\
	\end{array}
	\right):\ a\neq 0, c\in \mathbb{F} \right\},$
\end{itemize}
\end{theorem}
\textbf{Proof.}

Let $ A=A_{1}(\alpha_1, \alpha_2, \alpha_4, \beta_1)=\left(
\begin{array}{cccc}
 \alpha_1 & \alpha_2 &\alpha_2+1 & \alpha_4 \\
 \beta_1 & -\alpha_1 & -\alpha_1+1 & -\alpha_2
\end{array}\right).$\\
Due to (\ref{aut}) one has the system of equations:
\begin{equation} \label{SE1}
\begin{array}{ccccccc}
-a c+a \alpha _1-a^2 \alpha _1-2 a c \alpha _2-c^2 \alpha _4+b \beta _1 &=0,\\
-b c-b \alpha _1-a b \alpha _1+a \alpha _2-b c \alpha _2-a d \alpha _2-c
d \alpha _4 &=0, \\
a+b-a d-b \alpha _1-a b \alpha _1+a \alpha _2-b c \alpha _2-a d \alpha _2-c d \alpha _4 &=0,\\
-b d-b^2 \alpha _1-b \alpha _2-2 b d \alpha_2+a \alpha _4-d^2 \alpha _4 &=0,\\
-a c+c \alpha _1+2 a c \alpha _1+c^2 \alpha _2-a^2 \beta _1+d \beta _1 &=0,\\
-b c+b c \alpha _1-d \alpha _1+a d \alpha _1+c \alpha _2+c d \alpha _2-a b \beta _1 &=0,\\
c+d-a d+b c \alpha _1-d \alpha _1+a d \alpha _1+c \alpha _2+c d \alpha _2-a b \beta _1 &=0,\\
-b d+2 b d \alpha _1-d \alpha _2+d^2 \alpha_2+c \alpha _4-b^2 \beta _1&=0.\end{array}\end{equation}

The equations $2$ and $3,$ of the system (\ref{SE1}) imply $a+b-ad+bc=0,$ i.e., $b=\Delta-a.$
Similarly the equations $6$ and $7,$ of the system (\ref{SE1}) imply $c+d-ad+bc=0,$ i.e., $d=\Delta-c.$\\
Therefore, $\Delta=a(\Delta-c)-c(\Delta-a),$ and this implies
\[\Delta(1-a+c)=0.\] Since $\Delta\neq 0$ we have $c=a-1.$ The similar observation implies $d= b+1.$\\ Therefore we get $g=\left(\begin{array}{cc}a&b\\
	a-1&b+1\end{array}\right)$, with $\Delta=a+b\neq 0$.

As a result the system (\ref{SE1}) can be rewritten as follows:
\begin{equation} \label{SE2}
\begin{array}{cc}
a^2(1+ \alpha _1+2 \alpha _2+ \alpha _4)-a(1+ \alpha _1+2  \alpha _2+2 \alpha _4)+\alpha _4-b \beta _1 &=0,\\
ab(\alpha_1+2\alpha _2+\alpha_4+1)+a\alpha _4+b(\alpha_1-\alpha _2-\alpha_4-1)-\alpha _4 &=0,\\
b^2(1+ \alpha _1+2 \alpha _2+ \alpha _4)-\alpha_4a+(3\alpha _2+2\alpha _4+1)b+\alpha _4&=0, \\
a^2(1-2 \alpha _1- \alpha _2+ \beta _1)+a(-1+ \alpha _1+2  \alpha _2)-\alpha _2+\alpha_1-\beta _1(b+1) &=0,\\
ab(1-2\alpha _1- \alpha _2+ \beta _1)-a(\alpha _1+2 \alpha_2)+b(2\alpha _1+ \alpha _2 -1) +\alpha _1+2 \alpha _2 &=0,\\
b^2(1-2\alpha _1- \alpha _2+ \beta _1)+ b(-2\alpha _1-\alpha _2+1)-a \alpha _4+\alpha _4&=0.
\end{array}\end{equation}

Note that if $(a,b)$ is a solution to the system and $b=0$ then $a=0$ or $a=1$. Indeed, the substitution $0$ for $b$ in the system gives\\
\[
\begin{array}{cc}
a^2(1+ \alpha _1+2 \alpha _2+ \alpha _4)-a(1+ \alpha _1+2  \alpha _2+2 \alpha _4)+\alpha _4&=0,\\
\alpha _4(a-1) &=0,\\
a^2(1-2 \alpha _1- \alpha _2+ \beta _1)+a(-1+ \alpha _1+2  \alpha _2)-\alpha _2+\alpha_1-\beta _1 &=0,\\
(\alpha _1+2 \alpha_2)(a-1)&=0.
\end{array}\]
 Due to the first two equalities one has
$a(1+ \alpha _1+2 \alpha _2)(a-1)=0$, and therefore $a(a-1)=0$. In $a=0$ case $g$ is singular, and it is out of the consideration, but if $a=1$ then $g$ is the identity matrix. Therefore further it is assumed that $b(a+b)\neq 0$.

The equations $2$ and $5$ of the system of equations (\ref{SE2}) imply
\begin{equation}\label{4}(3 \alpha _1+3 \alpha _2 +\alpha _4-\beta_1)a b +( \alpha _1+2 \alpha _2+\alpha _4)a- ( \alpha _1+2 \alpha _2+\alpha _4)b-( \alpha _1+2 \alpha _2+\alpha _4) =0,\end{equation} the equations $3$ and $6$ imply $b(b(3 \alpha _1+3 \alpha _2 +\alpha _4-\beta_1) +2( \alpha _1+2 \alpha _2+\alpha _4)) =0,$ i.e.,
\begin{equation}\label{5}b(3 \alpha _1+3 \alpha _2 +\alpha _4-\beta_1) +2( \alpha _1+2 \alpha _2+\alpha _4) =0.\end{equation}
The equations $2$ and $3$ imply $b(b+a)(1+\alpha _1+2\alpha _2 +\alpha _4)+b(\alpha _1+2 \alpha _2+\alpha _4) =0,$ i.e.,
\begin{equation}\label{6}(b+a)(1+\alpha _1+2\alpha _2 +\alpha _4)+(\alpha _1+2 \alpha _2+\alpha _4)=b+a+(b+a+1)(\alpha _1+2\alpha _2 +\alpha _4)=0.\end{equation}

But (\ref{4}), (\ref{5}) imply \[(b+a+1)(\alpha _1+2\alpha _2 +\alpha _4)=0 \] and therefore due to (\ref{6}) one has $a+b=0$, which is a contradiction. Therefore
\begin{equation}\label{7}
Aut(A_1(\alpha_1,\alpha_2,\alpha_4,\beta_1))=\{I\}.\end{equation}

Note that the proof of (\ref{7}) does not depend on characteristics of $\mathbb{F}$.

Let $A=A_{2}(\alpha_1, \beta_1, \beta_2)=\left(
\begin{array}{cccc}
 \alpha_1 & 0 & 0 & 1 \\
 \beta _1& \beta _2& 1-\alpha_1&0
\end{array}\right).$\\
Due to (\ref{aut}) one has the system of equations:
\begin{equation} \label{SE3}
\begin{array}{ccccccc}
-c^2+a \alpha _1-a^2 \alpha _1+b \beta _1 &=0,\\
-c d-a b \alpha _1+b \beta _2 &=0,\\
b-c d-b \alpha _1-a b \alpha _1 &=0,\\
a-d^2-b^2 \alpha _1 &=0,\\
-a c+c \alpha _1+a c \alpha _1-a^2 \beta _1+d \beta _1-a c \beta _2 &=0,\\
-b c+b c \alpha _1-a b \beta _1+d \beta _2-a d \beta _2 &=0,\\
d-a d-d \alpha _1+ad \alpha _1-a b \beta _1-b c \beta _2 &=0,\\
c-b d+b d \alpha _1-b^2 \beta _1-b d \beta _2&=0,\end{array}\end{equation}
From the equations $2$ and $3$ of the system (\ref{SE3}) we get $b(1-\beta_2-\alpha_1)=0.$ The following cases occur:

\underline{\textbf{Case 1}. $1-\alpha_1-\beta_2\neq 0$.} Then $b=0$ and one gets $c=0$, so $a\alpha_1(1-a)=0$ and $d(1-a)(1-\alpha_1)=0$ due to the equations $1$ and $7,$ respectively. Thus $a=1.$ The equation $4$ implies $d^2=1$ and $\beta_1(d-a^2)=0$ due to the equation $5.$ We get two cases:

\ \ \underline{Case 1.1. $\beta_1=0.$} Then $d=\pm 1$, therefore $g$ is $\left(
\begin{array}{cc}
1 & 0 \\
0 & \pm 1 \\
\end{array}
\right).$

\ \ \underline{Case 1.2. $\beta_1\neq0.$} Then $d=1$, so $g=I.$

\underline{\textbf{Case 2}. $\beta_2=1-\alpha_1.$} If $b=0$ it is easy to see that $g$ equals $\left(
\begin{array}{cc}
1 & 0 \\
0 & \pm 1 \\
\end{array}
\right)$ if $\beta_1=0,$  and $g=I$ if $\beta_1\neq0.$ Now one can assume that $b\neq 0$, so due to the equation $4$ and $1$  we have $\alpha _1= \frac{a-d^2}{b^2},\ $ $\beta _1 = \frac{c^2 b^2-a^2+a d^2+a^3-a^2 d^2}{b^3},$ respectively. The substitution them into (\ref{SE3}) results the following system of equations:
\begin{equation} \label{SE4}
\begin{array}{ccccccc}
-a-a^2+b^2-b c d+d^2+a d^2 &=0,\\
a^4-a^5+a b c+2 a^2 b c-2 a b^3 c-a^2 b^2 c^2-a^2 d+a^3 d\\
+b^2 c^2 d-a^3 d^2+a^4 d^2-b c d^2-2 a b c d^2+a d^3-a^2 d^3 &=0,\\
a^3-a^4+a b c-b^3 c-a b^2 c^2-a d+a^2 d+b^2 d-a b^2 d\\
-a^2 d^2+a^3 d^2-b c d^2+d^3-a d^3 &=0,\\
a^2-a^3+b c-b^2 c^2+2 a d-2 b^2 d-a d^2+a^2 d^2-2 d^3 &=0.\end{array}\end{equation}
Now we make use the equations $2$ and $4$ to get
\[ a b c+a^2 b c-2 a b^3 c-a^2 d-a^3 d+2 a^2 b^2 d+b^2 c^2 d-b c d^2-2 a b c d^2+a d^3+a^2 d^3=0,\]
that is
\[b c(a+a^2+b c d- d^2-a d^2)-a d(a+a^2+b c d-d^2-a d^2)-2 a b^3 c+2 a^2 b^2 d=0.\]
Then the equation $1$ of the system gives
\[b^3 c-a db^2-2 a b^3 c+2 a^2 b^2 d=0\]
\[b^3 c(1-2a) -a db^2(1-2 a)=b^2(1-2a)(b c -a d)=0.\]
This implies $a=\frac{1}{2}.$ The substitution it into the system (\ref{SE4}) yields:
\begin{equation}\label{SE5}
\begin{array}{ccccccc}
-3+4 b^2-4 b c d+6 d^2 &= 0, \\
1+32 b c-32 b^3 c-8 b^2 c^2-4 d+32 b^2 c^2 d-2 d^2-64 b c d^2+8 d^3 &=0,\\
1+8 b c-16 b^3 c-8 b^2 c^2-4 d+8 b^2 d-2 d^2-16 b c d^2+8 d^3 &=0,\\
1+8 b c-8 b^2 c^2+8 d-16 b^2 d-2 d^2-16 d^3 &=0.
\end{array}\end{equation}
Due to the equations $1$ and $2$ of the system of equations (\ref{SE5}) one gets
\[
1+8 b c-8 b^2 c^2-4 d-2 d^2-18 b c d^2+8 d^3 =0
\]
and therefore taking into account the equation $3,$ one has
$8 b^2 c-4 b d =4 b(2 b c- d)=0,$ i.e., $d=2 b c.$ This implies
$ \Delta = a d- b c = b c - b c=0,$ i.e., $g$ is singular. Therefore \[Aut(A_2(\alpha_1,\beta_1,\beta_2))=\{I\}, \ \mbox{if} \ \beta_1\neq 0 ,\]
and
\[Aut(A_2(\alpha_1,0,\beta_2))=\left\{I, \left(
\begin{array}{cc}
1 & 0 \\
0 & -1 \\
\end{array}
\right)\right\}.\]\\
Let $ A=A_{3}(\beta_1, \beta_2)=\left(
\begin{array}{cccc}
 0 & 1 & 1 & 0 \\
 \beta _1& \beta _2 & 1&-1
\end{array}\right).$\\
Due to (\ref{aut}) one has the system of equations:
\begin{equation} \label{SE6}
\begin{array}{ccccccc}
-2 a c+b \beta _1 &=0,\\
a-b c-a d+b \beta _2 &=0,\\
a+b-b c-a d &=0,\\
-b-2 b d &=0, \\
-a c+c^2-a^2 \beta _1+d \beta _1-a c \beta _2 &=0,\\
c-b c+c d-a b \beta _1+d \beta _2-a d \beta _2 &=0,\\
c+d-a d+c d-a b \beta _1-b c \beta _2 &=0,\\
-d-b d+d^2-b^2\beta _1-b d \beta _2&=0.\end{array}\end{equation}
The equations $2$ and $3$ imply $b(1-\beta_2)=0.$ The following cases may occur:

\underline{\textbf{Case 1}. $\beta_2\neq 1.$ } In this case $b=0$, $c=0$ and $a(1-d)=0$, i.e.,  $d=1.$ Therefore $a=1$ and $g=I.$

\underline{\textbf{Case 2}. $\beta_2=1.$} If $b=0$ it is easy to see that $g=I.$ If $b\neq 0$ then the equation $4$ implies $d=-\frac{1}{2},$ hence (\ref{SE6}) becomes
\begin{equation} \label{SE7}
\begin{array}{ccccccc}
-2 a c+b \beta _1 &=0,\\
\frac{3 a}{2}+b-b c &=0,\\
-2 a c+c^2-\frac{\beta _1}{2}-a^2 \beta _1 &=0,\\
-\frac{1}{2}+\frac{a}{2}+\frac{c}{2}-b c-a b \beta _1 &=0,\\
\frac{3}{4}+b-b^2 \beta _1&=0.\end{array}\end{equation}

The equation $2$ of the system (\ref{SE7}) implies $a= \frac{2}{3}(-b+b c)$ and since $b\neq0$ one gets $\beta _1= \frac{2 a c}{b}.$ Therefore, the system of equations (\ref{SE7}) becomes
\[\begin{array}{ccccccc}
\frac{2 c}{3}+\frac{4 b c}{3}+\frac{16 b^2 c}{27}+\frac{c^2}{3}-\frac{4 b c^2}{3}-\frac{16 b^2 c^2}{9}+\frac{16 b^2 c^3}{9}-\frac{16 b^2 c^4}{27}&=0,\\
-\frac{1}{2}-\frac{b}{3}+\frac{c}{2}-\frac{2 b c}{3}-\frac{8 b^2 c}{9}+\frac{16 b^2 c^2}{9}-\frac{8 b^2 c^3}{9}&=0,\\
\frac{3}{4}+b+\frac{4 b^2 c}{3}-\frac{4 b^2 c^2}{3}&=0.\end{array}\]
The identity \[-\frac{1}{2}-\frac{b}{3}+\frac{c}{2}-\frac{2 b c}{3}-\frac{8 b^2 c}{9}+\frac{16 b^2 c^2}{9}-\frac{8 b^2 c^3}{9}=(\frac{3}{4}+b+\frac{4 b^2 c}{3}-\frac{4 b^2 c^2}{3})(\frac{2}{3}c-\frac{2}{3})+\frac{b(1-4c)}{3},\]
shows that $g$ may be only singular as far as in this case $\Delta=\frac{b(1-4c)}{3}$.

Let $ A=A_{4}(\alpha_1, \beta_2)=\left(
\begin{array}{cccc}
 \alpha _1 & 0 & 0 & 0 \\
 0 & \beta _2& 1-\alpha _1&0
\end{array}\right).$

Owing to (\ref{aut}) one has the system of equations:
\[\begin{array}{ccccccc}
a \alpha _1-a^2 \alpha _1 &=0,\\
-a b \alpha _1+b \beta _2 &=0,\\
b-b \alpha _1-a b \alpha _1 &=0,\\
-b^2 \alpha _1 &=0,\\
-a c+c \alpha _1+a c \alpha _1-a c \beta _2 &=0,\\
-b c+b c \alpha _1+d \beta _2-a d \beta _2 &=0,\\
d-a d-d \alpha _1+a d \alpha _1-b c \beta _2 &=0,\\
-b d+bd \alpha _1-b d \beta _2&=0\end{array}.\]
It is easy to see that for this system $b=0$ (so $a\neq 0$, $d\neq 0$ due to $\Delta\neq 0$),
 $a\alpha_1(1-a)=0,$ $d\beta_2(1-a)=0,$ and $d(1-a)(1-\alpha_1)=0$, hence $a=1.$ Therefore,
 $c(-1+2\alpha_1-\beta_2)=0$ by the equation $5.$

\underline{\textbf{Case 1}. $\beta_2=2\alpha_1-1.$} We get $g=\left(
\begin{array}{cc}
1 & 0 \\
c & d \\
\end{array}
\right),$ where $d\neq 0.$

\underline{\textbf{Case 2}. $\beta_2\neq2\alpha_1-1.$} In this case $c=0$ and we obtain $g=\left(
\begin{array}{cc}
1 & 0 \\
0 & d \\
\end{array}
\right),$ where $d\neq 0.$

Let $ A=A_{5}(\alpha_1)=\left(
\begin{array}{cccc}
 \alpha_1& 0 & 0 & 0 \\
 1 & 2\alpha_1-1 & 1-\alpha_1&0
\end{array}\right).$\\
 From (\ref{aut}) we have the system of equations:
\[\begin{array}{ccccccc}
b+a \alpha _1-a^2 \alpha _1 &=0,\\
-b+2 b \alpha _1-a b \alpha _1 &=0,\\
b-b \alpha _1-a b \alpha _1 &=0,\\
-b^2 \alpha _1&=0,\\
-a^2+d+c \alpha _1-a c \alpha _1 &=0,\\
-a b-b c-d+a d+b c \alpha _1+2 d \alpha _1-2 a d \alpha _1 &=0,\\
-a b+b c+d-a d-2 b c \alpha _1-d \alpha _1+a d \alpha_1 &=0,\\
-b^2-b d \alpha _1&=0\end{array}.\]
Here immediately we get $b=0$, $a=1$, $d=1$ and $g=\left(
\begin{array}{cc}
1 & 0 \\
c & 1 \\
\end{array}
\right).$

Let $ A_{6}(\alpha_1, \beta_1)=\left(
\begin{array}{cccc}
 \alpha_1 & 0 & 0 & 1 \\
 \beta _1& 1-\alpha_1 & -\alpha_1&0
\end{array}\right).$

Then (\ref{aut}) gives the system of equations:
\[\begin{array}{ccccccc}
-c^2+a \alpha _1-a^2 \alpha _1+b \beta _1 &=0,\\
b-c d-b \alpha _1-a b \alpha _1 &=0,\\
-c d-b \alpha _1-a b \alpha _1 &=0,\\
a-d^2-b^2 \alpha _1 &=0,\\
-a c+c \alpha _1+2 a c \alpha _1-a^2 \beta _1+d \beta _1 &=0, \\
d-a d+b c \alpha _1-d \alpha _1+a d \alpha _1-a b \beta _1 &=0,\\
-b c+b c \alpha _1-d \alpha_1+a d \alpha _1-a b \beta _1 &=0,\\
c-b d+2 b d \alpha _1-b^2 \beta _1&=0.\end{array}.\]
It is easy to see from the system that $c=0$, $b=0,$ and  $a=1.$ The equations $4$ and $5,\ $ imply $d^2=1$, $\beta_1(d-1)=0,$ therefore we have to consider the following two cases:

\underline{\textbf{Case 1}. $\beta_1=0.$} We get $d=\pm 1$ and $g=\left(
\begin{array}{cc}
1 & 0 \\
0 & \pm 1 \\
\end{array}
\right).$

\underline{\textbf{Case 2}. $\beta_1\neq 0.$} We obtain $d= 1$ and $g=I.$

Let $ A_{7}(\beta_1)=\left(
\begin{array}{cccc}
 0 & 1 & 1 & 0 \\
 \beta_1& 1& 0&-1
\end{array}\right).$  Then\\
$gA_{7}-A_{7}(g\otimes g)=$ \[\left(
\begin{array}{cccc}
-2 a c+b \beta _1 & a+b-b c-a d &
 a-b c-a d & -b-2 b d \\
-a c+c^2-a^2 \beta _1+d \beta _1 & c+d-a d+c d-a b \beta _1 & c-b c+c d-a b \beta _1 & -d-b d+d^2-b^2 \beta _1
\end{array}
\right),\] Thank to (\ref{aut}) we get $b=0,$ $d=1$, $c=0$ and $a=1,$ which implies $g=I.$

Let $ A=A_{8}(\alpha _1)=\left(
\begin{array}{cccc}
 \alpha_1 & 0 & 0 & 0 \\
 0 & 1-\alpha_1 & -\alpha_1&0
\end{array}\right).$ Then\\
$gA_{8}-A_{8}(g\otimes g)=$ \[\left(
\begin{array}{cccc}
a \alpha _1-a^2 \alpha _1 & b-b \alpha _1-a b \alpha _1 & -b \alpha _1-a b \alpha _1 & -b^2 \alpha _1 \\
-a c+c \alpha _1+2 a c \alpha _1 & d-a d+b c \alpha _1-d \alpha _1+a d \alpha _1 & -b c+b c \alpha _1-d \alpha _1+a d \alpha _1 & -b d+2 b d \alpha_1
\end{array}
\right),\] According to (\ref{aut}) one gets $b=0,$  $d(1-a)=0$, that is $a=1$, $-c+3c\alpha_1=0.$ Again we consider two cases:

\underline{\textbf{Case 1}. $\alpha_1=\frac{1}{3}.$} We get $g=\left(
\begin{array}{cc}
1 & 0 \\
c & d \\
\end{array}
\right),$ where $d\neq 0.$

\underline{\textbf{Case 2}. $\alpha_1\neq \frac{1}{3}.$} Then $c=0$ and  $g=\left(
\begin{array}{cc}
1 & 0 \\
0 & d \\
\end{array}
\right),$ where $d\neq0.$
Therefore,
\[ Aut(A_8(\alpha_1))=\left\{\left(
\begin{array}{cc}
1 & 0 \\
0 & d \\
\end{array}
\right):\ d\neq 0 \right\},\ \mbox{where}\  \alpha_1\neq \frac{1}{3},\]
 \[Aut\left(A_8\left(\frac{1}{3}\right)\right)=\left\{\left(
\begin{array}{cc}
1 & 0 \\
c & d \\
\end{array}
\right):\ d\neq 0 \right\}. \ \ \ \ \ \ \ \ \ \ \ \ \ \ \ \ \ \ \ \ \ \ \]

Let $ A=A_{9}=\left(
\begin{array}{cccc}
 \frac{1}{3}& 0 & 0 & 0 \\
 1 & \frac{2}{3} & -\frac{1}{3}&0
\end{array}\right).$ Then\\
\[gA_{9}-A_{9}(g\otimes g)=\left(
\begin{array}{cccc}
\frac{a}{3}-\frac{a^2}{3}+b & \frac{2 b}{3}-\frac{a b}{3} & -\frac{b}{3}-\frac{a b}{3} & -\frac{b^2}{3} \\
-a^2+\frac{c}{3}-\frac{a c}{3}+d & -a b+\frac{b c}{3}+\frac{2 d}{3}-\frac{2 a d}{3} & -a b-\frac{2 b c}{3}-\frac{d}{3}+\frac{a d}{3} & -b^2-\frac{b d}{3}
\end{array}
\right).\] Owing to (\ref{aut}) one has $b=0$, $a=1$,  $d=1$ and $g=\left(
\begin{array}{cc}
1 & 0 \\
c & 1 \\
\end{array}
\right).$

Let $ A=A_{10}=\left(
 \begin{array}{cccc}
 0 & 1 & 1 & 0 \\
 0 &0& 0 &-1
\end{array}
\right).$ Then\\
\[gA_{10}-A_{10}(g\otimes g)=\left(
\begin{array}{cccc}
-2 a c & a-b c-a d & a-b c-a d & -b-2 b d \\
c^2 & c+c d & c+c d & -d+d^2
\end{array}
\right).\] Due to (\ref{aut}) it is easy to see that $c=0$, $d=1$ and $b=0,$ i.e., $g=\left(
\begin{array}{cc}
a & 0 \\
0 & 1 \\
\end{array}
\right),$ where $a\neq 0$.

Let $ A=A_{11}=\left(
 \begin{array}{cccc}
 0 & 1 & 1 & 0 \\
 1 &0& 0 &-1
\end{array}
\right).$ Then\\
\[gA_{11}-A_{11}(g\otimes g)=\left(
\begin{array}{cccc}
b-2 a c & a-b c-a d & a-b c-a d & -b-2 b d \\
-a^2+c^2+d & -a b+c+c d & -a b+c+c d & -b^2-d+d^2
\end{array}
\right).\] The equation (\ref{aut}) gives $b(1+2d)=0.$ The following cases may occur:

\underline{\textbf{Case 1}. $b=0.$} Then one has $c=0,$ $d=1,$ $a=\pm 1$ and $g=\left(
\begin{array}{cc}
\pm 1 & 0 \\
0 & 1 \\
\end{array}
\right).$

\underline{\textbf{Case 2}. $d=-\frac{1}{2}.$} In this case   $b=\pm \frac{\sqrt{3}}{2}$, $a=\pm \frac{1}{2}$ and $c=\frac{b}{2a}$.

\ \ \underline{Case 2.1. $a=\frac{1}{2}.$} Then $c=\pm \frac{\sqrt{3}}{2}$ and $g=\left(
\begin{array}{cc}
\frac{1}{2} & \pm \frac{\sqrt{3}}{2} \\
\pm \frac{\sqrt{3}}{2} & -\frac{1}{2}  \\
\end{array}
\right).$

\ \ \underline{Case 2.2. $a=-\frac{1}{2}.$} Then $c=\mp \frac{\sqrt{3}}{2}$ and $g=\left(
\begin{array}{cc}
-\frac{1}{2} & \pm \frac{\sqrt{3}}{2} \\
\mp \frac{\sqrt{3}}{2} & -\frac{1}{2}  \\
\end{array}
\right).$

Let $ A=A_{12}=\left(
 \begin{array}{cccc}
 0 & 0 & 0 & 0 \\
 1 &0&0 &0\end{array}
\right).$ Then\\
\[gA_{12}-A_{12}(g\otimes g)=\left(
\begin{array}{cccc}
b & 0 & 0 & 0 \\
-a^2+d & -a b & -a b & -b^2
\end{array}
\right).\] Due to (\ref{aut}) one has
$g=\left(
\begin{array}{cc}
a & 0 \\
c & a^2 \\
\end{array}
\right)$, where $a\neq 0$.

Here are the corresponding results in the cases of characteristic $2$ and $3.$ The proof is similar to that of the case of characteristic not $2$ and $3.$

\begin{theorem}\label{thm5} The automorphism groups of the algebras listed in Theorem $\ref{thm2}$  are given as follows
	\begin{itemize}
\item $Aut(A_{1,2}(\alpha_1,\alpha_2,\alpha_4,\beta_1))=\{I\},$

\item $Aut(A_{2,2}(\alpha_1,\beta_1,\beta_2))=\{I\},$
\item $Aut(A_{3,2}(\alpha_1,\beta_2))=\left\{ I,\ \left(
\begin{array}{cc}
1 & 0 \\
1+\beta_2 & 1 \\
\end{array}
\right)\right\}$,
\item $Aut(A_{4,2}(\alpha_1,\beta_2))=\left\{ \left(
\begin{array}{cc}
1 & 0 \\
0 & d \\
\end{array}
\right):\ d\neq 0 \right\}$ if $\beta_2\neq 1 ,$
\item $Aut(A_{4,2}(\alpha_1,1 ))=\left\{ \left(
\begin{array}{cc}
1 & 0 \\
c & d \\
\end{array}
\right):\ c\in \mathbb{F}, d\neq 0 \right\},$\\
\item $Aut(A_{5,2}(\alpha_1))=\left\{ \left(
\begin{array}{cc}
1 & 0 \\
c & 1 \\
\end{array}
\right):\ c\in \mathbb{F} \right\},$\\
\item $Aut(A_{6,2}(\alpha_1,\beta_1))=\{I\},$
\item $Aut(A_{7,2}(\alpha_1)=\left\{ I,\ \left(
\begin{array}{cc}
1 &0 \\
1+\alpha_1 & 1 \\
\end{array}
\right)\right\},$
\item $Aut(A_{8,2}(\alpha_1))=\left\{\left(
	\begin{array}{cc}
	1 & 0 \\
	0 & d \\
	\end{array}
	\right):\ d\neq 0 \right\}$ if $\alpha_1\neq 1 ,$
 $Aut(A_{8,2}(1))=\left\{\left(
	\begin{array}{cc}
	1 & 0 \\
	c & d \\
	\end{array}
	\right):\ c\in \mathbb{F}, d\neq 0 \right\},$
\item $Aut(A_{9,2})=\left\{\left(
	\begin{array}{cc}
	1 & 0 \\
	c & 1 \\
	\end{array}
	\right):\ c\in \mathbb{F} \right\},$
\item $Aut(A_{10,2})=\left\{\left(
	\begin{array}{cc}
	a & 0 \\
	0 & 1 \\
	\end{array}
	\right):\ a\neq 0 \right\},$
\item $Aut(A_{11,2})=\left\{I,\ \left(
\begin{array}{cc}
0 & 1 \\
1 & 0 \\
\end{array}
\right),\ \left(
\begin{array}{cc}
0 &  1 \\
 1 & 1  \\
\end{array}
\right),\ \left(
\begin{array}{cc}
1 & 0 \\
1 & 1 \\
\end{array}
\right),\ \left(
\begin{array}{cc}
1 &  1 \\
1 & 0  \\
\end{array}
\right),\ \left(
\begin{array}{cc}
1 & 1 \\
 0 & 1 \\
\end{array}
\right)\right\},$
\item $Aut(A_{12,2})=\left\{\left(
	\begin{array}{cc}
	a & 0 \\
	c & a^2 \\
	\end{array}
	\right):\ a\neq 0, c\in \mathbb{F} \right\},$
\end{itemize}
\end{theorem}

\begin{theorem}\label{thm6} The automorphism groups of the algebras listed in Theorem $\ref{thm3}$  are given as follows
	\begin{itemize}
\item $Aut(A_{1,3}(\alpha_1,\alpha_2,\alpha_4,\beta_1))=\{I\},$

\item $Aut(A_{2,3}(\alpha_1,\beta_1,\beta_2))=\{I\}$ if $\ \beta_1\neq 0 ,$
 $Aut(A_{2,3}(\alpha_1,0,\beta_2))=\left\{I, \left(
\begin{array}{cc}
1 & 0 \\
0 & -1 \\
\end{array}
\right)\right\},$
\item $Aut(A_{3,3}(\beta_1,\beta_2))=\{I\},$
\item $Aut(A_{4,3}(\alpha_1,\beta_2))=\left\{ \left(
\begin{array}{cc}
1 & 0 \\
0 & d \\
\end{array}
\right):\ d\neq 0 \right\}$ if $\beta_2\neq 2\alpha_1-1 ,$
\item $Aut(A_{4,3}(\alpha_1,2\alpha_1-1 ))=\left\{ \left(
\begin{array}{cc}
1 & 0 \\
c & d \\
\end{array}
\right):\ c\in \mathbb{F}, d\neq 0 \right\},$\\
\item $Aut(A_{5,3}(\alpha_1))=\left\{ \left(
\begin{array}{cc}
1 & 0 \\
c & 1 \\
\end{array}
\right):\ c\in \mathbb{F}, \right\},$\\
\item $Aut(A_{6,3}(\alpha_1,\beta_1))=\{I\}$ if $\beta_1\neq 0 ,$
 $Aut(A_{6,3}(\alpha_1,0))=\left\{I, \left(
\begin{array}{cc}
1 & 0 \\
0 & -1 \\
\end{array}
\right)\right\},$
\item $Aut(A_{7,3}(\beta_1))=\{I\}.$
\item $Aut(A_{8,3}(\alpha_1))=\left\{\left(
	\begin{array}{cc}
	1 & 0 \\
	0 & d \\
	\end{array}
	\right):\ d\neq 0 \right\},$
\item $Aut(A_{9,3})=\left\{\left(
	\begin{array}{cc}
	1 & 0 \\
	c & 1 \\
	\end{array}
	\right):\ c\in \mathbb{F} \right\},$
\item $Aut(A_{10,3})=\left\{\left(
	\begin{array}{cc}
	a & 0 \\
	0 & 1 \\
	\end{array}
	\right):\ a\neq 0 \right\},$
\item $Aut(A_{11,3})=\left\{\ \left(
\begin{array}{cc}
 1 & 0 \\
c & 1 \\
\end{array}
\right):\ c\in \mathbb{F}\right\},$
\item $Aut(A_{12,3})=\left\{\left(
	\begin{array}{cc}
	a & 0 \\
	c & a^2 \\
	\end{array}
	\right):\ a\neq 0, c\in \mathbb{F} \right\},$
\end{itemize}
\end{theorem}
\begin{remark} Another interpretation of Theorem $\ref{thm4}$ \emph{(}Theorem $\ref{thm5},$ Theorem $\ref{thm6})$ is that the stabilizers, with respect to the action $(\ref{3}),$ of the matrices listed in Theorem $\ref{thm1}$ \emph{(}respectively, Theorem $\ref{thm2},$ Theorem $\ref{thm3})$ are described.\end{remark}

\section{ Derivations of $2$-dimensional algebras }

If $\mathbb{A}$ is an algebra given by MSC $A$ then, due to (\ref{2}) the algebra of its derivations $Der(A)$ is represented as follows
\begin{equation}\label{der}
  Der(A)=\{D\in M(2;\mathbb{F}):\  A(D\otimes I+I\otimes D)-DA=0\}.
\end{equation}
Further we use the notation $D=\left(
\begin{array}{cc}
	a & b \\
	c & d
\end{array}
\right).$

One of the main results of this section is given in the following theorem.
\begin{theorem}\label{thm7} The derivations of all algebra structures on $2$-dimensional vector space over an algebraically closed field  $\mathbb{F}$ of characteristic not $2,\ 3$ are given as follows.
\begin{itemize}
  \item $Der(A_{1}(\alpha_1, \alpha_2, \alpha_4, \beta_1))= Der(A_{2}(\alpha_1, \beta_1, \beta_2))= Der(A_{3}(\beta_1, \beta_2))=\{0\},$
  \item $Der(A_{4}(\alpha_1, \beta_2))=\left\{\left(
\begin{array}{cc}
  0 & 0 \\
  0 & d \\
 \end{array}
\right)
  :\ d\in \mathbb{F} \right\}$ if $ \beta_2\neq 2\alpha_1-1,$
  \item $Der(A_{4}(\alpha_1, 2\alpha_1-1))=\left\{\left(
  \begin{array}{cc}
    0 & 0 \\
    c & d \\
   \end{array}
  \right):\ c,\ d\in \mathbb{F}\right\},$
  \item $Der(A_{5}(\alpha_1))=\left\{\left(
\begin{array}{cc}
0 & 0 \\
c & 0 \\
\end{array} \right):\ c\in \mathbb{F}\right\},$
  \item $Der(A_{6}(\alpha_1, \beta_1))=Der(A_{7}(\beta_1))=\{0\},$
  \item $Der(A_{8}(\alpha_1))=\left\{\left(
\begin{array}{cc}
0 & 0 \\
0 & d \\
\end{array} \right):\ d\in \mathbb{F}\right\}$ if $\alpha_1\neq \frac{1}{3} ,$
 $Der(A_{8}(\frac{1}{3}))=\left\{\left(
\begin{array}{cc}
0 & 0 \\
c & d \\
\end{array} \right):\ c,\ d\in \mathbb{F}\right\},$
\item $ Der(A_{9})=\left\{\left(
\begin{array}{cc}
0 & 0 \\
c & 0 \\
\end{array} \right):\ c \in \mathbb{F}\right\} ,$
\item $Der(A_{10})=\left\{\left(
\begin{array}{cc}
a & 0 \\
0 & 0 \\
\end{array} \right):\ a \in \mathbb{F}\right\},$
\item $Der(A_{11})=\{0\},$
\item $Der(A_{12})=\left\{\left(
\begin{array}{cc}
a & 0 \\
c & 2a \\
\end{array} \right):\ a,\ c \in \mathbb{F}\right\} ,$
\end{itemize}
\end{theorem}
\textbf{Proof.}
Let $ A=A_{1}(\alpha_1, \alpha_2, \alpha_4, \beta_1)=\left(
\begin{array}{cccc}
 \alpha_1 & \alpha_2 &\alpha_2+1 & \alpha_4 \\
 \beta_1 & -\alpha_1 & -\alpha_1+1 & -\alpha_2
\end{array}\right).$\\
Then the equality
$ A_{1}(\alpha_1, \alpha_2, \alpha_4, \beta_1)(D\otimes I+I\otimes D)-DA_{1}(\alpha_1, \alpha_2, \alpha_4, \beta_1)=0$
is equivalent to the following system of equations
\begin{equation*} \label{SE8}
\begin{array}{cccccccc}
c+a \alpha _1+2 c \alpha _2-b \beta _1&=0,\\
 c-3 c \alpha _1+2 a \beta _1-d \beta _1&=0,\\
  2 b \alpha _1+d \alpha _2+c \alpha _4&=0,\\
   -a \alpha _1-2 c \alpha _2+b \beta _1&=0,\\
   -b+d+2 b \alpha _1+d \alpha _2+c \alpha _4&=0,\\
    a-c-a \alpha _1-2 c \alpha _2+b \beta _1&=0,\\
     b+3 b \alpha _2-a \alpha _4+2 d \alpha _4 &=0,\\
      b-2 b \alpha _1-d \alpha _2-c \alpha _4&=0\\
      \end{array}
      \end{equation*}
The equations $3$ and $8$ of the system of equations above imply $b=0,$ the equations $1$ and $6$ imply $a=0,$ the equations $3$ and $5$ imply $d=0$ and the equations $1$ and $4$ imply $c=0,$ therefore we get $D=0.$

Let $ A=A_{2}(\alpha_1, \beta_1, \beta_2)=\left(
\begin{array}{cccc}
 \alpha_1 & 0 & 0 & 1 \\
 \beta _1& \beta _2& 1-\alpha_1&0
\end{array}\right).$\\
Then the equation (\ref{der}) is equivalent to the system of equations:
\begin{equation*} \label{SE9}
\begin{array}{cccccccc}
a \alpha _1-b \beta _1&=0,\\
  c-2 c \alpha _1+2 a \beta _1-d \beta _1+c \beta _2&=0,\\
   c+b \alpha _1-b \beta _2&=0,\\
    b \beta _1+a \beta _2&=0,\\
 -b+c+2 b \alpha _1&=0,\\
  a-a \alpha _1+b \beta _1 &=0,\\
   -a+2 d &=0,\\
    b-c-b \alpha _1+b \beta _2&=0.
   \end{array}
   \end{equation*}
The equations $1$ and $6$ imply $a=0,$ then the equation $7$ gives $d=0,$ the equations $3$ and $8$ imply $b=0$ and then due to the equation $3$ one gets $c=0,$ therefore one has $D=0.$

Let $ A=A_{3}(\beta_1, \beta_2)=\left(
\begin{array}{cccc}
 0 & 1 & 1 & 0 \\
 \beta _1& \beta _2 & 1&-1
\end{array}\right).$\\
Then due to (\ref{der}) one has the system of equations:
\begin{equation*} \label{SE10}
\begin{array}{ccccccc}
 2 c-b \beta _1&=0,\\
  c+2 a \beta _1-d \beta _1+c \beta _2&=0,\\
   d-b \beta _2&=0,\\
    -2 c+b \beta _1+a \beta _2&=0\\
     -b+d&=0,\\
      a-2 c+b \beta _1&=0,\\
      3 b &=0,\\
        b-d+b \beta _2&=0.\end{array}\end{equation*}
The equations $5$ and $7$ imply $b=d=0$, the equation $1$ gives $c=0,$ the equation $6$ implies $a=0,$ hence $D=0.$

Let $ A=A_{4}(\alpha_1, \beta_2)=\left(
\begin{array}{cccc}
 \alpha _1 & 0 & 0 & 0 \\
 0 & \beta _2& 1-\alpha _1&0
\end{array}\right).$\\
Then the equation (\ref{der}) is equivalent to the system of equations:
\begin{equation*} \label{SE11}
\begin{array}{ccccccc}
a \alpha _1&=0,\\
c-2 c \alpha _1+c \beta _2&=0,\\
 b \alpha _1-b \beta _2&=0,\\
  a \beta _2&=0,\\
   -b+2 b \alpha _1 &=0,\\
    a-a \alpha _1&=0,\\
    b-b \alpha _1+b \beta _2&=0.\end{array}\end{equation*}
The equations $3$ and $7$ imply $b=0,$ the equations $1$ and $6$ imply $a=0.$ Therefore, $c(1-2\alpha_1+\beta_2)=0.$  If $\beta_2\neq 2\alpha_1-1$ one gets $D=\left(
                         \begin{array}{cc}
                           0 & 0 \\
                           0 & d \\
                         \end{array}
                       \right),$
if $\beta_2= 2\alpha_1-1$ we obtain $D=\left(
                         \begin{array}{cc}
                           0 & 0 \\
                           c & d \\
                         \end{array}
                       \right).$

Let $ A=A_{5}(\alpha_1)=\left(
\begin{array}{cccc}
 \alpha_1& 0 & 0 & 0 \\
 1 & 2\alpha_1-1 & 1-\alpha_1&0
\end{array}\right).$\\
Then due to (\ref{der}) one has the system of equations:
\begin{equation*} \label{SE12}
\begin{array}{ccccccc}
-b+a \alpha _1&=0,\\
 2 a-d&=0,\\
  b-b \alpha _1&=0,\\
   -a+b+2 a \alpha _1&=0,\\
    -b+2 b \alpha _1&=0,\\
a+b-a \alpha _1&=0,\\
 b \alpha _1&=0.\end{array}\end{equation*}
The equations $3$ and $7$ imply $b=0,$ the equations $1$ and $6$ give $a=0$, then according to the equation $2$ we have $d=0.$ Therefore, $D=\left(
                         \begin{array}{cc}
                           0 & 0 \\
                           c & 0 \\
                         \end{array}
                       \right).$

Let $ A_{6}(\alpha_1, \beta_1)=\left(
\begin{array}{cccc}
 \alpha_1 & 0 & 0 & 1 \\
 \beta _1& 1-\alpha_1 & -\alpha_1&0
\end{array}\right).$\\
Then due to (\ref{der}) one has the system of equations:
\begin{equation*} \label{SE13}
\begin{array}{ccccccc}
a \alpha _1-b \beta _1 &=0,\\
c-3 c \alpha _1+2 a \beta _1-d \beta _1&=0,\\
 -b+c+2 b \alpha _1 &=0,\\
  a-a \alpha _1+b \beta _1&=0,\\
 c+2 b \alpha _1 &=0,\\
  -a \alpha _1+b \beta _1&=0,\\
   -a+2 d &=0,\\
    b-c-2 b \alpha _1&=0.\end{array}\end{equation*}
The equations $1$ and $4$ imply $a=0,$ the equation $7$ yields $d=0,$ the equations $3$ and $5,$ imply $b=0,$ therefore $c=0.$ Hence, $D=0.$

Let $ A_{7}(\beta_1)=\left(
\begin{array}{cccc}
 0 & 1 & 1 & 0 \\
 \beta_1& 1& 0&-1
\end{array}\right).$
In this case the equation (\ref{der}) is equivalent to the system of equations:
\begin{equation*} \label{SE14}
\begin{array}{ccccccc}
2 c-b \beta _1&=0,\\
c+2 a \beta _1-d \beta _1&=0,\\
  -b+d &=0,\\
   a-2 c+b \beta _1&=0,\\
d&=0,\\
 -2 c+b \beta _1&=0,\\
  3 b&=0,\\
   b-d&=0.\end{array}\end{equation*}
We get $b=d=c=a=0$ and therefore, $D=0.$\\

Let $ A=A_{8}(\alpha _1)=\left(
\begin{array}{cccc}
 \alpha_1 & 0 & 0 & 0 \\
 0 & 1-\alpha_1 & -\alpha_1&0
\end{array}\right).$ Then
\[ A_8(\alpha _1)(D\otimes I+I\otimes D)-DA_8(\alpha _1)=\left(
\begin{array}{cccc}
 a \alpha _1 & -b+2 b \alpha _1 & 2 b \alpha _1 & 0 \\
 c-3 c \alpha _1 & a-a \alpha _1 & -a \alpha _1 & b-2 b \alpha _1
\end{array}
\right).\] We rewrite the equation (\ref{der}) in form of the system of equations as follows:
\begin{equation*} \label{SE15}
\begin{array}{ccccccc}
a \alpha _1&=0,\\
 c-3 c \alpha _1&=0,\\
  -b+2 b \alpha _1&=0,\\
   a-a \alpha _1&=0,\\
    2 b \alpha _1&=0,\\
     -a \alpha _1&=0,\\
      b-2 b \alpha _1&=0.\end{array}\end{equation*}
The equations $1$ and $4$ imply $a=0$, the equations $3$ and $5$ imply $b=0$ and then $c(1-3\alpha_1)=0.$ If $\alpha_1 \neq\frac{1}{3},$  one gets $D=\left(                                                         \begin{array}{cc} 0 & 0 \\
                  0 & d \\
                  \end{array}\right),$
if $\alpha_1 =\frac{1}{3},$  one has $D=\left(                                                                \begin{array}{cc}
                                                                  0 & 0 \\
                                                                  c & d \\
                                                                \end{array}
                                                              \right)$\\

Let $ A=A_{9}=\left(
\begin{array}{cccc}
 \frac{1}{3}& 0 & 0 & 0 \\
 1 & \frac{2}{3} & -\frac{1}{3}&0
\end{array}\right).$ Then
\[A_9(D\otimes I+I\otimes D)-DA_9=\left(
\begin{array}{cccc}
 \frac{a}{3}-b & -\frac{b}{3} & \frac{2 b}{3} & 0 \\
 2 a-d & \frac{2 a}{3}+b & -\frac{a}{3}+b & \frac{b}{3}
\end{array}
\right).\]
Due to (\ref{der}) one has $a=b=d=0$ and $D=\left( \begin{array}{cc}
                                                   0 & 0 \\
                                                   c & 0 \\
                                                 \end{array}
                                               \right)$\\

Let $ A=A_{10}=\left(
 \begin{array}{cccc}
 0 & 1 & 1 & 0 \\
 0 &0& 0 &-1
\end{array}
\right).$ Then \[ A_{10}(D\otimes I+I\otimes D)-DA_{10}=\left(
\begin{array}{cccc}
2 c & d & d & 3 b \\
0 & -2 c & -2 c & -d
\end{array}
\right).\] Due to (\ref{der}) one has $c=d=b=0$ and $D=\left(\begin{array}{cc} a & 0 \\
                                                        0 & 0 \\
                                                      \end{array}
                                                    \right).$

Let $ A=A_{11}=\left(
 \begin{array}{cccc}
 0 & 1 & 1 & 0 \\
 1 &0& 0 &-1
\end{array}
\right).$ Then
\[ A_{11}(D\otimes I+I\otimes D)-DA_{11}=\left(
\begin{array}{cccc}
 -b+2 c & d & d & 3 b \\
 2 a-d & b-2 c & b-2 c & -d
\end{array}
\right)=0.\] From the equation (\ref{der}) it is immediate that $D=0.$

Finally, for $ A=A_{12}=\left(
 \begin{array}{cccc}
 0 & 0 & 0 & 0 \\
 1 &0&0 &0\end{array}
\right).$ Then\\
$ A_{12}(D\otimes I+I\otimes D)-DA_{12}=\left(
\begin{array}{cccc}
 -b & 0 & 0 & 0 \\
 2 a-d & b & b & 0
\end{array}
\right)=0$, the equation (\ref{der}) gives
$D=\left(
  \begin{array}{cc}
      a & 0 \\
     c & 2a \\
      \end{array}
      \right).$

Now we present the corresponding results for characteristic $2$ and $3$ cases without any justifications as far as they are similar to that of the proofs in $Char(\mathbb{F})\neq 2,3$  case above.

\begin{theorem}\label{thm8} The derivations of all the algebras on $2$-dimensional vector space over an algebraically closed field  $\mathbb{F}$ of characteristic $2$ are given as follows.
\begin{itemize}
  \item $Der(A_{1,2}(\alpha_1, \alpha_2, \alpha_4, \beta_1))=\{0\},$
  \item $Der(A_{2,2}(\alpha_1, \beta_1, \beta_2))=\{0\}$ if $ \beta_1\neq 0,$
   $Der(A_{2,2}(\alpha_1, 0, \beta_2))= \left\{\left(
     \begin{array}{cc}
       0 & 0 \\
        0 & d \\
        \end{array}
        \right) :\ d\in \mathbb{F}\right\},$
  \item $Der(A_{3,2}(\alpha_1, \beta_2))=\left\{0\right\} $ if $\beta_2\neq 1,$
  $ Der(A_{3,2}(\alpha_1,1))=\left\{\left(
 \begin{array}{cc}
    0 & 0 \\
    c & 0 \\
     \end{array}
     \right):\ c\in \mathbb{F}\right\},$
  \item $Der(A_{4,2}(\alpha_1, \beta_2))=\left\{\left(
     \begin{array}{cc}
       0 & 0 \\
        0 & d \\
        \end{array}
        \right) :\ d\in \mathbb{F}\right\}$ if $\beta_2\neq 1,$
  \item $Der(A_{4,2}(\alpha_1, -1))=\left\{\left(
   \begin{array}{cc}
    0 & 0 \\
     c & d \\
     \end{array}
     \right):\ c,\ d\in \mathbb{F}\right\},$
  \item $Der(A_{5,2}(\alpha_1))=\left\{\left(
\begin{array}{cc}
0 & 0 \\
c & 0 \\
\end{array} \right):\ c\in \mathbb{F}\right\},$
  \item $Der(A_{6,2}(\alpha_1, \beta_1))=\{0\}$ if $ \beta_1\neq 0,$
  $Der(A_{6,2}(\alpha_1, 0))=\left\{\left(
 \begin{array}{cc}
            0 & 0 \\
            0 & d \\
           \end{array}\right):\ d\in \mathbb{F}\right\},$
 \item $Der(A_{7,2}(\alpha_1))=\{0\}$ if $\alpha_1\neq 1, $
 $Der(A_{7,2}(1))= \left\{\left( \begin{array}{cc}
  0 & 0 \\
 c & 0 \\
 \end{array}  \right):\ c\in \mathbb{F}\right\},$
 \item $ Der(A_{8,2}(\alpha_1))=\left\{\left(
\begin{array}{cc}
0 & 0 \\
0 & d \\
\end{array} \right):\ d\in \mathbb{F}\right\}$ if $ \alpha_1\neq 1 ,$
$Der(A_{8,2}(1))=\left\{\left(
\begin{array}{cc}
0 & 0 \\
c & d \\
\end{array} \right):\ c,\ d\in \mathbb{F}\right\},$
\item $Der(A_{9,2})=\left\{\left(
\begin{array}{cc}
0 & 0 \\
c & 0 \\
\end{array} \right):\ c \in \mathbb{F}\right\} ,$
$Der(A_{10,2})=\left\{\left(
\begin{array}{cc}
a & 0 \\
c & 0 \\
\end{array} \right):\ a,\ c \in \mathbb{F}\right\},$
\item $Der(A_{11,2})=\{0\},$
 $Der(A_{12,2})=\left\{\left(
\begin{array}{cc}
a & 0 \\
c & 0 \\
\end{array} \right):\ a,\ c \in \mathbb{F}\right\} ,$

\end{itemize}
\end{theorem}

\begin{theorem}\label{thm9} The derivation of all the algebras on $2$-dimensional vector space over an algebraically closed field  $\mathbb{F}$ of characteristic $3$ are given as follows.
\begin{itemize}
  \item $Der(A_{1,3}(\alpha_1, \alpha_2, \alpha_4, \beta_1))= Der(A_{2,3}(\alpha_1, \beta_1, \beta_2))= Der(A_{3,3}(\beta_1, \beta_2))=\{0\},$
  \item $ Der(A_{4,3}(\alpha_1, \beta_2))=\left\{\left(
\begin{array}{cc}
 0 & 0 \\
  0 & d \\
   \end{array}
    \right)
  :\ d\in \mathbb{F} \right\}$ if $ \beta_2\neq 2\alpha_1-1,$
  \item $Der(A_{4,3}(\alpha_1, 2\alpha_1-1))=\left\{\left(
 \begin{array}{cc}
  0 & 0 \\
c & d \\
\end{array}
\right):\ c,\ d\in \mathbb{F}\right\},$
$Der(A_{5,3}(\alpha_1))=\left\{\left(
\begin{array}{cc}
0 & 0 \\
c & 0 \\
\end{array} \right):\ c\in \mathbb{F}\right\},$
  \item $Der(A_{6,3}(\alpha_1, \beta_1))=Der(A_{7,3}(\beta_1))=\{0\},$
  $Der(A_{8,3}(\alpha_1))=\left\{\left(
\begin{array}{cc}
0 & 0 \\
0 & d \\
\end{array} \right):\ d\in \mathbb{F}\right\} ,$
  \item $Der(A_{9,3})=\left\{\left(
\begin{array}{cc}
0 & 2c \\
c & 0 \\
\end{array} \right):\ c \in \mathbb{F}\right\} ,$
 $Der(A_{10,3})=\left\{\left(
\begin{array}{cc}
a & b \\
0 & 0 \\
\end{array} \right):\ a,\ b \in \mathbb{F}\right\},$
  \item $Der(A_{11,3})=\left\{\left(
\begin{array}{cc}
0 & 0 \\
c & 0 \\
\end{array} \right):\ c \in \mathbb{F}\right\},$
  $Der(A_{12,3})=\left\{\left(
\begin{array}{cc}
a & 0 \\
c & 2a \\
\end{array} \right):\ a,\ c \in \mathbb{F}\right\} ,$
\end{itemize}
\end{theorem}
\begin{remark} The sets of the algebras of types $A_{1}(\alpha_1, \alpha_2, \alpha_4, \beta_1)$, $A_{1,2}(\alpha_1, \alpha_2, \alpha_4, \beta_1)$, $A_{1,3}(\alpha_1, \alpha_2, \alpha_4, \beta_1)$ are open, dense subsets of $V=M(m\times m^2;\mathbb{F})$ and therefore due to the results presented we can conclude that the majority of $2$-dimensional algebras have only trivial automorphisms and trivial derivations.	\end{remark}

\begin{center}{\textbf{Acknowledgments}}
\end{center}
The second authors research is supported by FRGS14-153-0394, MOHE and the third author acknowledges MOHE for supports by grant 01-02-14-1591FR.
\vskip 0.4 true cm


\begin{thebibliography}{99}
	
	\bibitem{A} H. Ahmed, U. Bekbaev, I. Rakhimov, Complete classification of two-dimensional algebras, AIP Conference Proceedings 1830, 070016 (2017); doi: 10.1063/1.4980965
		\bibitem{I} I. Kaygorodov, Y. Volkov, The variety of $2$-dimensional algebras over an algebraically closed field, \textit{ArXiv: 1701.08233[math.RA]}, 2017, 1-15.
\end{thebibliography}
\end{document}